\documentclass[epic,11pt]{article}
\usepackage{amsfonts,amsmath,amscd,amssymb,amsthm}
\usepackage{fullpage}
\usepackage[alphabetic,initials]{amsrefs}
\usepackage[all]{xy}
\usepackage{hyperref} 
\usepackage{mathrsfs}
\input xy%
\usepackage{fullpage}
\usepackage{amssymb}

\def\convf{\hbox{\space \raise-2mm\hbox{$\textstyle      \bigotimes \atop \scriptstyle \omega$} \space}}
\def\0{{\bar 0}}
\def\1{{\bar 1}}

\def\Z{{\mathbb Z}}
\def\F{{\mathbb F}}
\def\Q{{\mathbb Q}}
\def\B{{\mathcal B}}

\def\Vc{{V^c}}

\def\Max{{\operatorname{Max}\;}}
\def\Spec  {{\operatorname{Spec  }\;}}

\def\rad{\operatorname{rad\;}}
\def\Ker {{\operatorname{Ker}\;}}

\def\sign{{\operatorname{sign}}}
\newcommand{\ttk}{\mathtt{k}}

\newcommand{\ttT}{\mathtt{T}}

\newcommand{\itema}{\item[{{\rm$($a$)$}}]}
\newcommand{\itemb}{\item[{{\rm$($b$)$}}]}
\newcommand{\itemc}{\item[{{\rm$($c$)$}}]}
\newcommand{\itemd}{\item[{{\rm$($d$)$}}]}

\newcommand{\noi}{\noindent}
\newcommand{\ga}{\alpha}

\newcommand{\gc}{\gamma}

\newcommand{\Gc}{\Gamma}
\newcommand{\Gl}{\Lambda}

\newcommand{\gd}{\delta}

\newcommand{\gs}{\sigma}

\newcommand{\gl}{\lambda}

\newcommand{\gep}{\epsilon}

\newcommand{\op}{\oplus}

\newcommand{\A}{\mathcal A}

\newcommand{\ot}{\otimes}

\newcommand{\fg}{\mathfrak{g}}\newcommand{\fgl}{\mathfrak{gl}}
\newcommand{\fpsl}{\mathfrak{psl}}\newcommand{\osp}{\mathfrak{osp}}

\newcommand{\fh}{\mathfrak{h}}

\newcommand{\fb}{\mathfrak{b}}

\newcommand{\fa}{\mathfrak{a}}

\newcommand{\ff}{\footnote}
\newfont{\eufm}{eufm10 scaled\magstep1}
\newcommand{\pd}{\partial}
\newcommand{\ci}{\circ} \newcommand{\ti}{\times}

\newcommand{\cI}{\mathcal{I}}

\newcommand{\cB}{\mathcal{B}}

\newcommand{\cA}{\mathcal{A}}

\newcommand{\cS}{\mathcal{S}}

\newcommand{\cV}{\mathcal{V}}

\newcommand{\cW}{\mathfrak{W}}

\newcommand{\ey}{\end{eqnarray}}
\newcommand{\by}{\begin{eqnarray}}
\newcommand{\nn}{\nonumber}

\newcommand{\bco}{\begin{conjecture}}
\newcommand{\ba}{\begin{alg}}
\newcommand{\ea}{\end{alg}}
\newcommand{\eco}{\end{conjecture}}
\newcommand{\bpf}{\begin{proof}}
\newcommand{\epf}{\end{proof}}
\newcommand{\bt}{\begin{theorem}}

\newcommand{\et}{\end{theorem}}

\newcommand{\br}{\begin{rem}}
\newcommand{\er}{\end{rem}}
\newcommand{\brs}{\begin{rems}}
\newcommand{\ers}{\end{rems}}
\newcommand{\bi}{\begin{itemize}}
\newcommand{\ei}{\end{itemize}}
\newcommand{\bl}{\begin{lemma}}
\newcommand{\bsul}{\begin{sublemma}}
\newcommand{\esul}{\end{sublemma}}
\newcommand{\bp}{\begin{proposition}}
\newcommand{\be}{\begin{equation}}
\newcommand{\bc}{\begin{corollary}}
\newcommand{\bexs}{\begin{examples}}
\newcommand{\eexs}{\end{examples}}
\newcommand{\bexa}{\begin{example}}
\newcommand{\eexa}{\end{example}}
\newcommand{\bex}{\begin{exercise}}
\newcommand{\eex}{\end{exercise}}
\newcommand{\btab}{\begin{tab}}
\newcommand{\etab}{\end{tab}}
\newcommand{\el}{\end{lemma}}
\newcommand{\ep}{\end{proposition}}
\newcommand{\ee}{\end{equation}}
\newcommand{\ec}{\end{corollary}}
\newcommand{\Bc}{\begin{center}}
\newcommand{\Ec}{\end{center}}
\newcommand{\bh}{\begin{hyp}}
\newcommand{\eh}{\end{hyp}}
\newcommand{\bhs}{\begin{hyps}}
\newcommand{\ehs}{\end{hyps}}
\newcommand{\bd}{\begin{dfn}}
\newcommand{\ed}{\end{dfn}}

\begin{document}
\title{Table of Contents}

\newtheorem{thm}{Theorem}[section]
\newtheorem{hyp}[thm]{Hypothesis}
 \newtheorem{hyps}[thm]{Hypotheses}
  \newtheorem{rems}[thm]{Remarks}

\newtheorem{conjecture}[thm]{Conjecture}
\newtheorem{theorem}[thm]{Theorem}
\newtheorem{theorem a}[thm]{Theorem A}
\newtheorem{example}[thm]{Example}
\newtheorem{examples}[thm]{Examples}
\newtheorem{corollary}[thm]{Corollary}
\newtheorem{rem}[thm]{Remark}
\newtheorem{lemma}[thm]{Lemma}
\newtheorem{sublemma}[thm]{Sublemma}
\newtheorem{cor}[thm]{Corollary}
\newtheorem{proposition}[thm]{Proposition}
\newtheorem{exs}[thm]{Examples}
\newtheorem{ex}[thm]{Example}
\newtheorem{exercise}[thm]{Exercise}
\numberwithin{equation}{section}%
\setcounter{part}{0}
\newcommand{\drar}{\rightarrow}
\newcommand{\lra}{\longrightarrow}
\newcommand{\rra}{\longleftarrow}
\newcommand{\dra}{\Rightarrow}
\newcommand{\dla}{\Leftarrow}
\newcommand{\rl}{\longleftrightarrow}

\newtheorem{Thm}{Main Theorem}


\newtheorem*{thm*}{Theorem}
\newtheorem{lem}[thm]{Lemma}
\newtheorem*{lem*}{Lemma}
\newtheorem*{prop*}{Proposition}
\newtheorem*{cor*}{Corollary}
\newtheorem{dfn}[thm]{Definition}
\newtheorem*{defn*}{Definition}
\newtheorem{notadefn}[thm]{Notation and Definition}
\newtheorem*{notadefn*}{Notation and Definition}
\newtheorem{nota}[thm]{Notation}
\newtheorem*{nota*}{Notation}
\newtheorem{note}[thm]{Remark}
\newtheorem*{note*}{Remark}
\newtheorem*{notes*}{Remarks}
\newtheorem{hypo}[thm]{Hypothesis}
\newtheorem*{ex*}{Example}
\newtheorem{prob}[thm]{Problems}
\newtheorem{conj}[thm]{Conjecture}
\title{The Nullstellensatz for  supersymmetric polynomials.}
\author{Ian M. Musson\ff{Research partly supported by  Simons Foundation grant 318264.} \\Department of Mathematical Sciences\\
University of Wisconsin-Milwaukee\\ email: {\tt
musson@uwm.edu}}
\maketitle
\begin{abstract} In this paper we prove a Nullstellensatz for  supersymmetric polynomials.  This gives a bijection between radical ideals and superalgebraic sets.  These are algebraic sets which are invariant under the Weyl groupoid of Sergeev and Veselov, \cite{SV2}. Note that  the algebra of supersymmetric polynomials is not Noetherian, so the usual Nullstellensatz does not apply. 
However it deos satisfy   the ascending chain condition on radical ideals and this allows for the decomposition of superalgebraic sets into irreducible components.  
Analogous results hold for the a ring of Laurent supersymmetric polynomials. 
\\ \\
As   an application, we give a proof of conjecture 13.5.1 from \cite{M}.  This  concerns  the maximal ideals in the enveloping algebra of the general linear and  orthosymplectic  Lie superalgebras. The center is closely related to 
the algebra of supersymmetric polynomials  and the result can be thought of as an analog of the weak Nullstellensatz.  
\end{abstract}
\section{Introduction} \label{s10.3}
\noi  The interest in supersymmetric polynomials comes from several  sources.  First they often satisfy analogs of combinatorial 
properties of symmetric polynomials.  For example there is an analog of the Jacobi-Trudi identity for Schur polynomials \cite{PT}, and 
super-Schur polynomials form a $\Z$-basis for the $\Z$-algebra of supersymmetric polynomials. Secondly the center $Z(\fg)$ of the enveloping algebra $U(\fg)$ is isomorphic to an algebra of supersymmetric polynomials, when $\fg$ is the Lie superalgebra
 $\fgl(m,n)$ or the orthosymplectic Lie superalgebra $\osp(2m+1,2n)$. For the
 Lie superalgebra $\osp(2m,2n)$ a slight modification of this statement is necessary. Last{}ly it is shown in \cite{SV2} that if $\fg=\fgl(m|n)$, then a natural quotient of the Grothendieck group of finite dimensional $\fg$-modules is isomorphic to a ring of Laurent supersymmetric polynomials. Suprisingly though, the rings of (Laurent) supersymmetric polynomials have received little interest from a geometric point of view, perhaps because they are not Noetherian.  Yet as we show the spectra of these rings  have some very pleasant geometric properties. In fact they are quotient spaces, except that the quotient arises from the action of a groupoid rather than a group.

\section{General results.}
We begin with some general results. These will be applied to both supersymmetric polynomials and Laurent supersymmetric polynomials.  Unless otherwise stated all rings are commutative. 
\noi  If $I$ is an ideal of $A$, let $$\rad(I) = \{f\in A| f^n \in I, \mbox{ for some } n>0\}.$$
 If $I=\rad(I)$ we say $I$ is a {\it radical ideal}.

\bl  \label{owl} Suppose $A$ is a subring of $B$, 
$T\in B$ and $TB \subseteq A$. Let $P$ be an ideal 
of $A$  with $T\in P$ and suppose $P$ is a radical ideal, 
then $TB \subseteq P$
\el
\bpf This holds since $(TB)^2 =T^2B \subseteq AT \subseteq A P = P$.  
   \epf
\bexa The result may be false if $P$ is not a radical ideal.  Let $$A=\ttk +T\ttk[S,T],  \quad B =\ttk[S, T].$$
If $P = \ttk T +T^2\ttk[S,T]$, then 
$T\in P$ but 
 $TB$ is not contained in $P$.\eexa
We keep the assumptions of Lemma \ref{owl} and suppose that $T$ is a non-zero divisor in $B$.   Since $TB  \subset A \subset B$, we have $A_T =B_T $. 
If $L$ is an ideal of $A$ define the {\it  extension} of $L$ to $B_T$ to be $L^e = L_T$, the localization of $L$  with respect to the powers of $T$.  If $M$ is an ideal of $A_T$, set $M^c= M \cap A$, the {\it  contraction} of  $M$ to $A$.  By 
\cite{GlWd} Theorem 10.15, contraction and extension provide inverse lattice isomorphisms between
the lattice ideals of $ A_T$  and the lattice of ideals 
C of $A$ such that $A/C$ is $T$-torsionfree. Another result we need is that the contraction of a prime ideal is prime.  Indeed if $R$ is a subring of the commutative ring  $S$ and $P$ is prime in $R$, then $R/R\cap P$ embeds in $S/P$ which is a domain, so $R\cap P$ is prime. We mention this rather obvious fact  only because it fails for non-commutative non-Noetherian rings \cite{GlWd} Exercise 10.M. It is also easy to  check that if $P$ is prime in $R$, then $P_T$ is prime in $R_T$.  
Thus  extension  and contraction  give a bijection 

\be \label{yak} \{P \in \Spec   A| T \notin P\} \leftrightarrow \{ \Spec   A_T \}.\ee 
\noi Note that if $P$ is maximal, then by exactnesss of localization  \cite{R} Corollary 4.81, $P_T$  is also maximal.   From now on we identify the two sets on either side of \eqref{yak}.

\bc  \label{pig}  If         $m$ is a maximal ideal of  $A$ amd $T\notin m$, the  maximal ideal $M$ of  $B$ 
given by $M= m_T\cap B$ satisfies  $m = M \cap A$.    Also $mB \neq B$.     \ec \bpf    This first statement follows from Equation  \eqref{yak}, and the second is an immediate consequence.         \epf

\bt \label{cow} Let $\phi: A\lra A/BT= C$ be the natural map.  Then we have a disjoint union 
\be \label{gnu} \Spec   A = \Spec B_T \cup \phi^{-1} (\Spec   C).  \ee
\et 
\bpf As noted above we have $A_T =B_T $.    For a prime ideal $P$ of $A$ there are two possibilities.  If $T\notin P$, then $P_T$ is a prime ideal of $A_T=B_T$ such that $P_T\cap A=P.$ If $T\in P$, then $P$ is the inverse image of the  prime ideal of $P/BT$   under $\phi.$
\epf

  \brs{\rm\bi\itema  Clearly  $ \Spec   B_T$ is open and $\phi^{-1} (\Spec   C)$ is closed in $ \Spec   A$. In our applications  $B$ will be  a finitely generated algebra over $\Z$ or over a field, so $B_T$ will be a  Jacobson ring, \cite{E} Theorem 4.19, and then  $\Spec  B_T$  will be locally closed.
\itemb Theorem \ref{cow} actually holds for noncommutative rings provided $T$ is a {\it normal element}, that is $TB=BT$,  which is not a zero divisor. 
\ei}\ers
\noi Equation \eqref{gnu} does not hold if $\Spec$ is replaced everywhere by Rad, the set of radical ideals of  a ring. Just consider the intersection $P\cap Q$ of two primes lying in different components on the  right side of Equation \eqref{gnu} to see this.  
To understand the relationship between prime and radical ideals, it would help to know that every radical ideal is a finite intersection
of prime ideals.  By \cite {K} Theorem 87, this holds if the ring satisfies the ascending chain condition on radical ideals. We abbreviate this condition to ACCR. 
We are now ready for  the following result.

\bl \label{hos} In the situation of Theorem \ref{cow}, suppose that the rings 
$B_T $ and $ C$ satisfy ACCR.  Then so does $A$. \el
\bpf    Let $R^{1} \subseteq R^{2} \subseteq \ldots $ 
be an ascending chain of radical ideals in $A$.   
Then $R_T^{1} \subseteq R_T^{2}  \ldots $ is an ascending chain of radical ideals in $A_T=B_T$, so by assumption, there is an $m$ such that $R_T^{m} = R_T^{i}$ for all $i\ge m$.  
Write $R_T^{m} = p_1 \cap p_2 \cap  \ldots  \cap p_r$ for some prime ideals of $B_T$.
  We can assume that this intersection is irredundant, and all the $p_k$ are minimal over $R_T^{m} $.
Then if  $P_k= p_k \cap A$ and $i\ge m$,  the $P_k$ are exactly the minimal primes over $R^{i} $
 that do not contain $T$. By \cite{E} Corollary 2.12 every radical ideal in a commutative ring is the intersection of a (possibly infinite number) of prime ideals. So for $i\ge m$, write $R^{i}  $ as an intersection of prime ideals in $A$.  Say
\[ R^{i}  = P_1 \cap \ldots \cap P_r \cap  \bigcap_{j\in \Lambda_i} Q_j.
\]
 where the $Q_j$ are prime ideals that contain $TB$. Thus  $D_i =  \bigcap_{j\in \Lambda_i} Q_j$ is a radical ideal in $A$ containing $TB$, and if $\bar{D_i} = D_i/TB$, then $\bar{D}_{1} \subseteq \bar{D}_{2} \subseteq \ldots $ is an ascending chain of radical ideals of $C$.  Therefore there is an $n \ge m$, such that $\bar{D}_{i} = \bar{D}_{n}  $ for all $i>n$.  Hence $R^{i} = R^{n}  $ for all  $i>n$.  
\epf

\section{Supersymmetric Polynomials.}

We use the notation of \cite{M}. Fix nonnegative integers $m,n$. 
Let $\mathcal{X}_m = (x_1, \ldots, x_m),$ $ \mathcal{Y}_n = (y_1, \ldots, y_n)$ be two
 sets of indeterminates and  $L$  a commutative ring. 
   The symmetric group $\mathcal{S}_m$ acts on the
 polynomial ring $L [\mathcal{X}_m] = L [x_1, \ldots, x_m]$ by the rule
 \[ w(x_i) = x_{w(i)} . \]
 Similarly $W_{m,n} = \mathcal{S}_m \times \mathcal{S}_n$ acts on $L [\mathcal{X}_m, \mathcal{Y}_n] =
 L [x_1, \ldots, x_m, y_1, \ldots, y_n]$.
 For $f \in L [\mathcal{X}_m, \mathcal{Y}_n]$ and $t \in
 L,$ we write $f(x_1 = t, y_1 = -t)$ for the polynomial obtained by substituting
  $x_1 = -y_1 = t$  in $f.$\\  \\
A $W_{m,n}$-invariant polynomial $f$ is
{\it supersymmetric} if $f(x_1 = t, y_1 = -t)$ is a polynomial which
is independent of $t$.   
Equivalently by \cite{M} Lemma 12.1.1,
\be \label{ono} \frac{\pd f}{\pd x_1} - \frac{\pd f}{\pd y_1} \in (x_1 + y_1). \ee
For example if $r \geq 1$ the power sum
$p^{(r)}_{m,n}$ defined by
 \be \label{sow} p^{(r)}_{m,n} = (x^r_1 + \ldots + x^r_m) + (-1)^{r-1}(y^r_1 +
 \ldots + y^r_n) \ee is supersymmetric.
The set of all supersymmetric
 polynomials in with coefficients in $L$ is called the {\it $L$-algebra of
 supersymmetric polynomials in} $\mathcal{X}_m,\mathcal{Y}_n$.
 We denote this algebra by  $I_L(x_1, \ldots,   x_m| y_1, \ldots, y_n)$ or $I_L[\mathcal{X}_m, \mathcal{Y}_n]$.\\
 \\
\noi The ring $I_\Z[\mathcal{X}_m, \mathcal{Y}_n]$ has a $\Z$-basis consisting of certain super Schur polynomials.   Chapter 12 of \cite{M} gives three equivalent definitons of these polynomials.  For our purposes it is most conventient to use 
the definition of super Schur polynomials via symmetrization operators. We denote these polynomilas by $ F_\lambda(\mathcal{X}_m/\mathcal{Y}_n)$.  In order to define them we need some combinatorics. 
\\ \\
\noi 
 First, denote the Vandermonde determinant  by 
\[ \Delta(\mathcal{X}_m) = \prod_{1 \leq i < j \leq m} (x_i - x_j) .\]

\noi The {\it Young  diagram} of the partition $\lambda$ is the set of
points
\[ D_\lambda = \{(i,j) \in \mathbb{N} ^2| 1 \leq j \leq \lambda_i \}
.\]
  We represent this diagram by a set of boxes in the fourth
  quadrant, with the first index $i$ corresponding to rows and the
  second $j$ to columns.  
 Next the {\it $(m,n)$-hook} is the set  \[\{(i,j) \in \mathbb{N} ^2| i \leq m \;
\mbox{or}\; j \leq n \} ,\] and we say that a partition $\lambda$
is {\it contained in the $(m,n)$-hook} if
$ D_\lambda$ is contained in this subset,
 or equivalently if $\lambda_{m+1} \leq n$. 
\noi Denote by ${\mathcal H}(m,n)$ the set
 of partitions $\lambda$ contained in the $(m,n)$-hook. Let $D^{m,n}_\lambda$ be the following subset of $D_\lambda$
 \[  D^{m,n}_\lambda = \{ (i,j)| i \leq m, j \leq n, j \leq
 \lambda_i\} .\]
 Suppose $\lambda$ is contained in the $(m,n)$-hook.  The part of
 $D_\lambda$ outside $D^{m,n}_\lambda$ is determined by two
 partitions $\mu, \nu$ defined by
 \[
\mu_i = \max\{0, \lambda_i -
 n\},\quad \nu_j = \max\{0, \lambda'_j - m \},
\]
 where $\gl'$ is the transpose of the partition $\gl$. 
   The part of
 $D_\lambda$ to the right of the line $j = n$ is a translate of
 the diagram $D_\mu$, while the part of $D_\lambda$ below the line
 $i = m$ is a translate of $D_{\nu'}$.
\\ \\
Now define
 \[ g_\lambda(\mathcal{X}_m,\mathcal{Y}_n) = \prod^m_{i=1} x^{\mu_i + m - i}_i
 \prod^n_{i=1} y^{\nu_i + n - i}_i \prod_{(i,j)\in D^{m,n}_\lambda}
 (x_i + y_j),\]
  and
  \be \label{cub}   F_\lambda (\mathcal{X}_m/\mathcal{Y}_n) = \sum_{w \in W_{m,n}} w \left[
  \frac{g_\lambda(\mathcal{X}_m,\mathcal{Y}_n)}{\Delta(\mathcal{X}_m)\Delta(\mathcal{Y}_n)}
  \right] . \ee
  If $\lambda$ is not contained in the $(m,n)$-hook, set
  $F_\lambda(\mathcal{X}_m/\mathcal{Y}_n) = 0$. 
\\ \\
 Assume $m > 0$ and $n > 0$ and set \[{\mathcal H}^0_{m,n} = {\mathcal H}(m,n)\backslash ({\mathcal H}(m,n-1) \cup {\mathcal H}(m-1,n)).
\ff{There is an erratum in this definition on page 248 of \cite{M}.  For a  complete list of current known errata see  https://www.ams.org/publications/authors/books/postpub/gsm-131}
\]
\noi Note that a partition $\lambda$ in ${\mathcal H}(m,n)$ belongs to
 ${\mathcal H}^0_{m,n}$ if and only if $(m,n) \in D_\lambda$.  Define
  \be \label{ant} T=T(\mathcal{X}_m, \mathcal{Y}_n) = \prod^m_{i=1} \prod^n_{j=1}(x_i + y_j).\ee If 
$\gl \in {\mathcal H}^0_{m,n}$, we have 

 \[ g_\lambda(\mathcal{X}_m,\mathcal{Y}_n) = T \prod^m_{i=1} x^{\mu_i + m - i}_i
 \prod^n_{i=1} y^{\nu_i + n - i}_i ,\] and since $T$ is $W$-invariant, \eqref{cub} becomes

  \by \label{bug}   F_\lambda (\mathcal{X}_m/\mathcal{Y}_n) &=& T \sum_{w \in W_{m,n}} w \left[
  \frac{ \prod^m_{i=1} x^{\mu_i + m - i}_i
 \prod^n_{i=1} y^{\nu_i + n - i}_i }{\Delta(\mathcal{X}_m)\Delta(\mathcal{Y}_n)}
  \right]\nn \\ &=& T S_\mu (\mathcal{X}_m)  S_\nu (\mathcal{Y}_n), \ey
where $S_\mu,   S_\nu$ are the usual Schur polynomials. \\ \\
We summarize the main properties of the super Schur polynomials $F_\lambda (\mathcal{X}_m/\mathcal{Y}_n) $.

\bp \label{dog} Assume $\lambda$ is contained in the $(m,n)$-hook. Then we have 
\bi 
\itema The polynomials $F_\lambda(\mathcal{X}_m/\mathcal{Y}_n)$ are
 supersymmetric.
\itemb The $\mathbb{Z} $-algebra of supersymmetric polynomials in
$\mathcal{X}_m$
  and $\mathcal{Y}_n$ has a $\mathbb{Z} $-basis consisting of the
  $F_\lambda(\mathcal{X}_m/\mathcal{Y}_n)$ as $\lambda$ ranges over partitions
  contained in the $(m,n)$-hook.
 \itemc  The polynomials $F_\lambda(\mathcal{X}_m/\mathcal{Y}_n)$ with $\lambda \in {\mathcal
 H}^0_{m,n}$ form a basis for the $\mathbb{Z} $-module of
 supersymmetric polynomials in $\mathcal{X}_m, \mathcal{Y}_n$ for which the
 substitution $x_m = y_n = 0$ yields the zero polynomial
\itemd
Specializing $y_{n+1} = 0$ in the polynomial
$F_\lambda(\mathcal{X}_m/\mathcal{Y}_{n+1})$, or $x_{m+1} = 0$ in
the polynomial $F_\lambda(\mathcal{X}_{m+1}/\mathcal{Y}_n)$ yields
the polynomial $F_\lambda(\mathcal{X}_m/\mathcal{Y}_n)$.
 \ei\ep
\noi \bpf See \cite{M} Lemmas 12.2.5, 12.2.6 and Proposition 12.2.7. \epf

\bc  \label{cat}  The homomophsim 
\be \label{egg} \phi:I_\Z[\mathcal{X}_m, \mathcal{Y}_n] \lra I_\Z[\mathcal{X}_{m-1}, \mathcal{Y}_{n-1}] \ee
defined by the substitution ${x}_m =y_n =0$ 
 is surjective and has kernel equal to $T\Z[\mathcal{X}_m,\mathcal{Y}_n]^{W_{m,n}}$. \ec\bpf
By (d)
$\phi(F_\lambda(\mathcal{X}_m/\mathcal{Y}_n))= F_\lambda(\mathcal{X}_{m-1}/\mathcal{Y}_{n-1})$, if $\gl$ is in the $(m-1,n-1)$ hook. It follows from (b) that $\phi$ is surjective. The second statement follows from (c) and \eqref{bug}. \epf


\section{Supersymmetric Laurent polynomials and Grothendieck rings of basic classical Lie superalgebras.}%
We recall some results of Sergeev, \cite{Se},  changing some of the notation slightly for convenience. The changes are noted in the footnotes.  First define  the  algebra of {\it  Laurent symmetric polynomials} to be
$$
\Lambda_{m}=\Bbb Z[x_{1}^{\pm1},\dots,x_{m}^{\pm1}]^{\cS_{m}}.
$$
 Let  $\lambda_1,\dots,\lambda_m$ be a non-increasing sequence of integers. We define the  Euler character  $E_{\lambda}(x)$ by means of the following formula
$$
E_{\lambda}(x)\Delta_m(x)=\{x_1^{\lambda_1+m-1}\dots x_m^{\lambda_m}\},
$$
where  $\Delta_m(x)=\prod_{i<j}(x_i-x_j)$ and the brackets $\{\}$ mean alternation over the  group  $\cS_m$,
$$
\{f\}=\sum_{\sigma\in \cS_m} \sign (\sigma)\sigma(f)).
$$
The following result is \cite{Se} Theorem 2.3.
\begin{thm}\label{jay}  The $E_{\lambda}(x)$ with $\lambda_1,\dots,\lambda_m$  a non-increasing sequence of integers form a $\Z$-basis of the ring  $\Lambda_{m}$.
\end{thm}
\noi Next the ring  \ff{In \cite{Se} $\Lambda_{m,n}$ is denoted $\Lambda^{\pm}_{m,n}$, but we want to add further superscripts later.  A similar remark applies to the ring $\Lambda_m$ defined above}
$$
\Lambda_{m,n}=\{f\in \Bbb Z[x_1^{\pm1},\dots,x_m^{\pm1},y_1^{\pm1},\dots,y_n^{\pm1}]^{\cS_m\ti \cS_n}\mid x_i\frac{\partial f}{\partial x_i}+y_j\frac{\partial f}{\partial y_j}\in (x_i-y_j)\}
$$
 will be called the {\it  ring of Laurent supersymmetric polynomials.}
Instead of \eqref{ono} we have
\[  x_1\frac{\pd f}{\pd x_1} + y_1\frac{\pd f}{\pd y_1} \in (x_1 - y_1). \]
This accounts for the sign differences in the definitions of the maps \eqref{egg} and \eqref{elk}.

\noi If $mn>0$ then by Lemma \ref{cod} in the Appendix, we have the homomorphism
 \be \label{elk} \phi: \Lambda_{m,n}\longrightarrow \Lambda_{m-1,n-1},
\ee
defined by setting $x_m$  equal to $y_n$. 
Denote by $Q(m, n)$ the set of pairs of sequences of non-increasing  integers  $(\lambda,\mu) \in \Bbb Z^m \ti \Bbb Z^n$.
For $(\lambda,\mu) \in Q(m, n)$ set 
\ff{In \cite{Se} there is a more general definition valid for pairs of sequences of non-increasing  integers  $(\lambda,\mu)$ which are  not in $\Bbb Z^m \ti \Bbb Z^n$.}
\begin{equation}\label{eel}
 K_{\lambda,\mu}=
\prod_{j=1}^n\prod_{i=1}^m\left(1-\frac{y_{j}}{x_{i}}\right)E_{\lambda}(x_{1},\dots, x_{m})E_{\mu}(y_1\dots,y_n), 
\end{equation}
A $\Z$-basis for 
$\Lambda_{m,n}$ is given in \cite{Se}, 
Theorem 5.6. From the proof  we have the following important statements.
\begin{thm}\label{fox}  The map $\phi$ is onto, and the $K_{\lambda,\mu}\,$
with $(\lambda,\mu)\in Q(m,n)$ 
form a $\Z$-basis for  $\Ker \phi$. 
\end{thm}
\noi Now set 
\be \label{rat}T= K_{0,0}=\prod_{j=1}^n\prod_{i=1}^m\left(1-\frac{y_{j}}{x_{i}}\right),\ee
$\A =\Lambda_{m,n}$, $\A'= \Lambda_{m-1,n-1}$ and $\B = \Bbb Z[x_1^{\pm1},\dots,x_m^{\pm1},y_1^{\pm1},\dots,y_n^{\pm1}]^{\cS_m\ti \cS_n} = \Lambda_m \ti \Lambda_n.
$
Then from Theorem \ref{jay} and \eqref{eel}, $T\B \subset \A  \subset \B$ and hence $\A_T=\B_T$. Also by Theorem \ref{fox},
 $\Ker \phi= T\B$. Therefore by  Equation  \eqref{gnu} we have a disjoint union

\[\Spec  \A = \Spec  \B_T \cup \phi^{-1} (\Spec  \A').\]
Part of the interest in the ring $\Lambda_{m,n}$ comes from the following result of Sergeev and Veselov, \cite{SV2}.  
Let $K(\mathfrak{g})$  be the quotient of the Grothendieck ring of finite dimensional $\Bbb Z_2-$graded  representations of the Lie superalgebra 
$\mathfrak{g} =  \mathfrak{gl}(m|n)$  by the ideal generated by all  $[M]+ [\Pi M]$ where $\Pi$ is the parity change functor, and 
$[M]$ is the class of the module $M$ in the Grothendieck ring.   
Then the supercharacter yields an isomorphism from $K(\mathfrak{g})$ to $\Lambda_{m,n}$.
  One of the key techniques of \cite{Se}
is the  map $\phi: \Lambda_{m,n} \rightarrow \Lambda_{m-1,n-1}$  from \eqref{elk}. It was shown by Hoyt and Reif, \cite{HR}  that this evaluation homomorphism has a natural interpretation  using the  Duflo–Serganova functor. This map as well as its anolog \eqref{egg} for supersymmetric
 polynomials is also one of the
main technical tools of this paper. 
\section{Prime  Ideals in Rings of Supersymmetric Polynomials.} \label{ai}
Now we apply our results simultaneously to rings of supersymmetric polynomials and rings of supersymmetric Laurent  polynomials. However some definitions still need to be made separately.  For this purpose we refer to case S or case L when dealing with 
supersymmetric polynomials or supersymmetric Laurent  polynomials respectively.  Set $W=W_{m,n}$ and $W'=W_{m-1,n-1}$. Then in case S, let 
$\A=\A_{m,n}$ be the ring of supersymmetric polynomials, and $\B= \Z[\mathcal{X}_m,\mathcal{Y}_n]^W.$ 
If $m\le 0$ or $n\le 0$, set $\A=\Z$.
 In addition set $ \A'={\A_{m-1,n-1}}$ and $\B'= \Z[\mathcal{X}_{m-1},\mathcal{Y}_{n-1}]^{W'}.$ The map $\phi$ is defined in Equation \eqref{egg} and the element $T$ is as in Equation \eqref{ant}.
\\ \\
In case L, we set 
$\A= \Lambda_{m,n}$ be the ring of supersymmetric Laurent polynomials, and $\B=\Lambda_{m} \ti \Lambda_{n}.$ 
If $m\le 0$ or $n\le 0$, set $\A=\A_{m,n}=\Z$.
 In addition set $ \A'={ \Lambda_{m-1,n-1}}$ and $\B'= \Lambda_{m-1} \ti \Lambda_{n-1}.$ The map $\phi$ is defined in Equation \eqref{elk} and the element $T$ in Equation \eqref{rat}.
In both cases we have a commutative diagram.

\begin{xy} 
(0,20)*+{}="f";
(40,20)*+{\A}="a";
 (100,20)*+{\A'}="b"; (40,0)*+{\B}="c"; 
(100,0)*+{\B'}="d"; {\ar "a";"b"}?*!/_2mm/{\phi}; {\ar "a";"c"};
{\ar "b";"d"};
{\ar@{>} "c";"d"};?*!/_2mm/{\Phi};
 \end{xy}


\noi The map $\Phi$ is given by evaluation of a polynomial at $x_m=y_n =0$. The two vertical maps are the obvious inclusions.

\br {\rm We make a remark about extension of scalars.  For  $\Z$-algebras $C, K$, set $C^K = C\ot_\Z K$.  If $K$ is a field of characteristic zero, then 
$I_\Z[\mathcal{X}_m, \mathcal{Y}_n]^K = I_K[\mathcal{X}_m, \mathcal{Y}_n].$ 
This is shown by first extending scalars to $\Q$ by localization, then to $K$ using a $\Q$-basis for the field extension $\Q\subseteq K$. Furthermore we can extend scalars in the above diagram  replacing each algebra $C$ by $C^K$, and the maps $\phi, \Phi$ by $\phi_K =\phi\ot 1,  \Phi_K=\Phi \ot 1.$ This also works when $K=L[T]$ where $L$ is a field of characteristic zero.  However it is not clear how to lift a supersymmetric polynomial with coefficients in $\F_p$ to the integers.  The polynomial $x^2+y^2$ is not supersymmetric, but mod 2 it is equal to $p^{(2)}_{1,1}$. Similar remarks apply to extension of scalars for $ \Lambda_{m,n}$.} \er

\noi For the remainder of this section $K$ will denote either the ring of integers, or  a field of characteristic zero.

\bc  If $\phi_K:\A^K =\A^K_{m,n}\lra \A^K_{m-1,n-1}$ is the map induced by Equation \eqref{egg} or, Equation \eqref{elk},
we have a disjoint union
\be \label{kit} \Spec  \A^K_{m,n} = \Spec  \B^K\cup \phi_K^{-1} (\Spec  \A^K_{m-1,n-1}).  \ee
\ec
\bpf Note that $T\B^K  \subset \A^K \subset \B^K$. Also we have $\Ker \phi_K = T\B^K$. If $K= \Z$, this follows from Corollary \ref{cat}. For a field argue as in the remark.    Hence the result follows from Equation  \eqref{gnu}.
\epf
\noi  To prove the Nullstellensatz we need more information on maximal ideals. 

\bp \label{doe} If  $m\in \Max \A^K_{m,n}$, then $m\B^K$  is a proper ideal of $\B^K$. \ep
\bpf If $T\notin m$ this follows from Corollary \ref{pig}.  If  $T\in m$, we use the above  commutative diagram.  Here $\phi_K(m)$ is  a maximal ideal of $\A^K_{m-1,n-1}$, hence  by induction  $J=\phi_K(m)\B^K_{m-1,n-1}$ is a proper ideal of $\B^K_{m-1,n-1}$. Since $\Phi_K(m\B^K)=J,$ it follows that $m\B^K$ 
is a proper ideal of  $\B^K$.
\epf
\bt\label{hen}
  If $m$ is a maximal ideal of  $\A^K$, then there is a maximal ideal $m$ of  $\B^K$ such that $m = M \cap \A^K$.\et  
\bpf If $T\notin m$ this follows from Corollary \ref{pig} extending scalars if necessary.  If $T\in m$, the result holds by taking $M$ to be any maximal ideal containing $m\B^K$.
\epf
\br {\rm The ideal $M$ is not unique, but in case S, if $\ttk$  is  an algebraically closed field the fibers of  the map $$\Max \B^\ttk_{m,n} \lra \Max \A^\ttk_{m,n}$$ given by $M\lra  M \cap \A^\ttk_{m,n}$  are known, \cite{M} Theorem 13.5.4. }\er
\bt \label{cob}
  In both cases S and L, the algebras $\A^K_{m,n} $  satisfy ACCR for $K =\Z$  or any field of characteristic zero. Hence any radical ideal $I$ is a finite intersection of prime ideals. These prime ideals can be taken to be the prime ideals minimal over $I$.\et
\bpf This follows from Lemma \ref{hos} extending scalars as necessary  (compare also Equation \eqref{kit}), once we observe that $B^K_T $  is Noetherian so has ACCR, and that $\A^K_{m-1,n-1} $ has ACCR by induction. \epf

\section{The Strong Nullstellensatz.}

\noi First we prove the analog of the weak Nullstellensatz.
Let $\ttk$ be an algebraically closed field of characteristic zero. Then 
set  $\A= \A^\ttk_{m,n}, \B= \B^\ttk_{m,n}$ 
and $ \ttk^{m|n} = \ttk^m\ti \ttk^n$
in case S.   
In case L, we set 
$\A= \Lambda^\ttk_{m,n}$,  $\B=\Lambda^\ttk_{m} \ti \Lambda^\ttk_{n}$  and $\ttk^{*m|n} =( \ttk^*)^m\ti ( \ttk^*)^n$. 
\bt If  $m$ is a maximal ideal of $\A$, there  exists $\gl \in \ttk^{m|n}$
in case S, or $\gl \in \ttk^{*m|n}$ in case L 
 such that $m = \{ f\in \cA| f(\gl)=0\}.$

 \et
\bpf This follows from Theorem \ref{hen}.   Note that  maximal ideals in $\B$ correspond to $W$-orbits in $\ttk^{m|n}$ or $\ttk^{*m|n}$. \epf
\noi If $I$ is a subset of $\A$, let
$\cV(I) = \{x\in \ttk ^{m|n}| f(x) =0 \mbox{ for all } f \in I\}$ in  case S, and $\cV(I) = \{x\in \ttk ^{*m|n}| f(x) =0 \mbox{ for all } f \in I\},$ in case L.  Such a set is called an {\it superalgebraic set}. If instead 
$I$ is a subset of $\B$, we say that 
$\cV(I)$ is an{\it algebraic set}.  Thus 
any superalgebraic set is algebraic.  
In addition if 
$V$ is a subset of $\ttk ^{m|n}$ in case S, or $\ttk^{*m|n}$ 
in case L, 
 set $$\cI_\cA(V) = \{f\in \cA| f(x) =0 \mbox{ for all } x \in V\}.$$  
We will also need 
$$\cI(V) = \{f\in \cB| f(x) =0 \mbox{ for all } x \in V\}.$$  
\bt \label{boa}  The maps $\cI_\cA$ and $\cV$ are inverse bijections between the set of superalgebraic sets in $\ttk ^{m|n}$  $($case S$)$ or $\ttk ^{*m|n}$  $($case L$)$, and the set of radical ideals in $\A$. Both maps are order reversing, \et
\bpf The key point is that $\cI_\cA(\cV(I)) \subseteq \rad(I)$. We adapt  a well-known argument of Rabinowitsch see  \cite{F} Chapter 1. 
Suppose $G\in \cI_\cA(\cV(I))$, and let $J$ be the ideal of $\A_{m,n} \ot \ttk [T]$ generated by $I$ and $TG-1$. Then $\cV(J)$ is empty, since $G$ vanishes whenever all polynomials in $I$ vanish. Therefore by the weak Nullstellensatz, $1\in J$, and we can write
\[1 = \sum_{i=1}^r A_i F_i + B(TG-1)\] where the $F_i $ are in $I$ and $B, A_i \in \cA_{m,n} \ot \ttk [T].$  
The result follows by multiplying by a large power of $G$. 
\epf
\bp \label{nit}
The maps $\cV,$ and $\cI_\cA$ satisfy  the following properties. Suppose that $E_\gl, V_\gl$ are subsets of  $\A$, and 
$\ttk ^{m|n}$ $($ or $\ttk ^{*m|n})$ 
respectively and that 
$\fa, \fb$ are ideals of $\cA.$ 
\bi
\itema $\cV(\cup_{\gl\in \Gl} E_\gl) = \cap_{\gl\in \Gl}\cV( E_\gl))$. 
\itemb $\cV(\fa \cap \fb) = \cV(\fa \fb) = \cV(\fa) \cup \cV(\fb)$.
\itemc $\cI_\cA(\cup_{\gl\in \Gl}V_\gl) = \cap_{\gl\in \Gl}\cI_\cA( V_\gl)$.
\ei\ep
\bpf Left to the reader. \epf

\section{ Superalgebraic sets.}
To give more meaning to Theorem \ref{boa} we need to know what  superalgebraic sets look like.  We show they are just the algebraic sets which are invariant under the Weyl groupoid of Sergeev and Veselov, \cite{SV2}. First we review their work.

\subsection{Definition of the Weyl Groupoid.}
In \cite{SV2} Sergeev and Veselov associated a certain groupoid $\mathfrak{W}= \mathfrak{W}(R)$,  which they the call Weyl groupoid,  to any  generalized root system 
$R\subset V$ in the  sense of Serganova \cite{S}.
A groupoid $\mathfrak G$  can be
defined as a small category with all morphisms  invertible.  
We denote the set of objects by $\mathfrak B$ which we call the {\it base}. 
As in \cite{SV2}  we use the same notation $\mathfrak G$ for the set of morphisms  as for the groupoid itself.
\\ \\
To define the Weyl groupoid we need a preliminary construction, namely
the semi-direct product groupoid $\Gc \ltimes \mathfrak G$.   Let $\mathfrak G$ be a groupoid and  $\Gamma$ a group  acting on $\mathfrak G$  by automorphisms of the corresponding category. In particular, $\Gamma$ acts on the base $\mathfrak B$ of $\mathfrak G$. 
Then the  {\it semi-direct product groupoid}
$\Gamma \ltimes \mathfrak G$ has the same base $\mathfrak B$, and the morphisms from $x$ to $y$  are pairs
$(\gamma,\,f),$ with $ \gamma \in \Gamma, f \in \mathfrak G$ such that $f : \gamma  x\rightarrow y.$
Note that if $(\gd,\,g)$ is a second morphism with $g : \gd y\rightarrow z$, then we also have morphisms $\gd(f) : \gd  \gamma x\rightarrow \gd y,$ and 
$g\ci\gd(f) : \gd  \gamma x \rightarrow z.$  
Hence we can define composition of morphisms  as follows: 

\[(\gd,\,g) \circ (\gc,\,f) = (\gd \gc,\,  g\ci\gd(f)).\]
Now we can  define the Weyl groupoid $\mathfrak{W}(R)$ corresponding to the generalized root system $R$.  Recall that the reflections with respect to the non-isotropic roots generate a finite group denoted $W_0.$ 
First  consider the following groupoid $\mathfrak T_{iso}$ with base $R_{iso},$ the set of all the isotropic roots in $R.$ The set of morphisms  $\alpha \rightarrow \beta$ is non-empty if and only if $\beta = \pm \alpha$ in which case it consists of just one element. Denote the corresponding morphism $\alpha \rightarrow -\alpha$ by $\tau_{\alpha}, \alpha \in R_{iso}.$
The group $W_0$ acts on $\mathfrak T_{iso}$ in a natural way: $\alpha \rightarrow w(\alpha),\,
\tau_{\alpha} \rightarrow \tau_{w(\alpha)}.$ 
Define the {\it Weyl groupoid} $$\mathfrak{W}(R) = W_0 \coprod  W_0 \ltimes \mathfrak T_{iso}$$ to be the disjoint union of the group $W_0$ considered as a groupoid with a single point base $[W_0]$ and the semi-direct product groupoid $W_0 \ltimes \mathfrak T_{iso}$ with base $R_{iso}.$  Observe that the disjoint union is a well defined operation on the groupoids.

 \subsection{Action of $\mathfrak{W}(R)$ on the ambient space $V$.}
For any set $X$ consider  the following groupoid $\mathfrak S(X)$, with base consisting  of all possible subsets $Y \subset X$ and morphisms are all possible bijections between them. By an {\it action of a groupoid $\mathfrak G$ on a set} $X$ we will mean a natural transformation between the categories $\mathfrak G$  and $\mathfrak S(X)$. If $X$ is a vector space, the   {\it affine groupoid} $\mathfrak A(X)$ has base all  affine subspaces, and  morphisms all affine bijections. Then an {\it affine action of $\mathfrak G$  on $X$} is 
natural transformation from $\mathfrak G$  to $\mathfrak A(X)$.
\\ \\
Returning to our  generalized root system 
$R\subset V$,  let $X= V$ and define the following affine action of the Weyl groupoid
$\mathfrak{W}(R)$ on $V$. 
The base point $[W_0]$ maps to the whole space $V$, 
meaning that an element $w \in W_0$ acts on any point of $V$ in the natural way. 
The base element corresponding to an isotropic root $\alpha$ maps to the hyperplane $\Pi_{\alpha}$ defined by the equation $(\alpha,x)=0.$
The element $\tau_{\alpha}$ acts as the shift 
$$\tau_{\alpha}(x) = x + \alpha, \, x \in \Pi_{\alpha}.$$
Note that since $\alpha$ is isotropic, $x + \alpha$ also belongs to $\Pi_{\alpha}.$ 
\\ \\
Let $V = \mathfrak h^*$ be the dual space to a Cartan subalgebra $\mathfrak h$
of a basic classical Lie superalgebra $\mathfrak g$ with generalized root system $R.$
 In the case of $\fpsl(n|n)$ we consider $\fgl(n|n)$ instead. 
Using the invariant bilinear form we can identify $V$ and $V^* = \mathfrak h$
and consider the elements of the group ring $\mathbb{Z}[\mathfrak {h}^*]$ as functions on $V.$
Similarly the symmetric algebra $S(\fh)$ consists of functions on $\fh^*$. 
A function $f$ on $V$ is {\it invariant under the action of the groupoid} $\mathfrak{W}$
if for any $g \in \mathfrak{W}$ we have $f(g(x)) = f(x)$ for all $x$ in  the domain of definition  of the morphism $g.$
\subsection{Description of Superalgebraic sets.}

\noi When dealing with a prime ideal $P$ of $\cA$ or $\cB$ with $T \notin P$, it is helpful to introduce  
$\ttk^{m,n} \ti \ttk$ or $\ttk^{*m,n} \ti \ttk$ using $z$ as the coordinate for the last copy of $\ttk$.  Then for $R= \cA$ or $\cB$ we have $R_T \cong R[z]/(Tz -1)$. Under this isomorphism $P_T$ corresponds to $P[z]/(Tz -1)$. 
A point $y \in
\ttk^{m,n}\mbox{ with } T(y) \neq 0$ 
 \ff{from this point on substitute $\ttk^{*m,n}$ for $\ttk^{m,n}$    if necessary.} corresponds to the point $(y, T(y)^{-1}) \in \ttk^{m,n} \ti \ttk$ .
This allows us to identify 

\[\cV(P_T) =\{ y \in \ttk^{m,n} | f(y) = 0 \mbox{ for all } f \in P, T(y)\neq 0 \}\] with 
\[\{ y \in \ttk^{m,n}\ti \ttk | f(y) = 0 \mbox{ for all } f \in P[z], (Tz-1)(y)\neq 0 \}.\] 
Note that $T$ is fixed by every element of $W=W_{m,n}$, so defining $wz=z$ for all $w \in W$ makes the map 
\be \label{air} y \lra (y, T(y)^{-1})\ee 
when $T(y) \neq 0,$ $W$-equivariant.  Also $T(\Pi_\ga)= 0$ for all $\ga \in R_{iso},$  the action of 
$W_0 \ltimes \mathfrak T_{iso}$ on
\[\cV_T =\{ y \in \ttk^{m,n} | T(y)\neq 0 \}\]  is trivial (that is no element of $\cV_T$ is in the domain of definition of any morphism in $W_0 \ltimes \mathfrak T_{iso}$).  Thus defining the action of $W_0 \ltimes \mathfrak T_{iso}$ on $\ttk^{m,n} \ti \ttk$ to be trivial, actually makes the map in \eqref{air}
$\mathfrak{W}(R)$-equivariant.

\bt\label{ass} The superalgebraic sets are exactly the algebraic sets that are invariant under the Weyl groupoid  $\mathfrak{W}$.\et
\bpf      There are two things to check.  If $I$ is an ideal of $\A^\ttk_{m,n}$ then $\cV(I)$  is invariant under  $\mathfrak{W}$.  This holds because elements of $I$ are supersymmetric, so $\cV(I)$ is a union of orbits. 
Conversely, suppose that $V$ is an algebraic set in the usual sense, which is invariant under $\cW$.   We want to show that $V$ is superalgebraic.  Use the notation before Theorem \ref{boa}. Set 
\[    V'= \{ v\in V| T(v) \neq 0\},\mbox{ and } \Vc:=V \backslash V' \]
\bi \itema Suppose that $P_1, \ldots , P_r$ are the prime ideals of $\cB$ minimal over $\cI(V)$ such that $T \notin P_i$. Then set $Q_i = P_i \cap \cA$, a prime ideal of $\cA$.  Using the fact that $T \notin Q_i$, it is easy to see that $P_i = Q_i\cB$.  Therefore $\cV(P_i) = \cV(Q_i).$ Using the remarks preceding the proof this shows that $V'$ is superalgebraic.  

\itemb 
 Now set 
$J = \cI(V) +T\cB$. 
Then 
 $\cV(J) =\Vc.$ This corresponds to an algebraic sets in $\Spec \cB'$. Moreover $\Vc$ is $\cW$-invariant so is superalgebraic by induction.  

\ei
It follows from Proposition \ref{nit} that $V$ is superalgebraic.
\epf
\subsection{Prime ideals and irreducible components.}
We say that a superalgebraic set is {\it irreducible} 
 if it cannot be written as the union of two proper superalgebraic subsets.

\bp  \label{hog}  If $I$ is a radical ideal and $V = \cV(I)$ is the corresponding superalgebraic set, then $I$ is prime if and only if $V$ is irreducible.\ep
\bpf  The proof is completley analagous to the classical case \cite{F} Proposition 1, page 7.\epf
\noi 
In general if 
$I$ is a radical ideal of $\cA$, then using Theorem \ref{cob}, we can write $I$ uniquely in the form $I =  P_1 \cap \ldots \cap P_r$,  
where the $P_i$ are the prime ideals of $\cA$ which are minimal over $I$.  In this situation we call the superalgebraic  sets $\cV(P_i)$ the  {\it irreducible components of $\cV(I)$.} By Proposition \ref{hog} they are irreducible.

\bc Every superalgebraic set is uniquely a finite union of  irreducible components. 
\ec

\bpf This follows  from  the Nullstellensatz and Proposition \ref{nit} (b).  
\epf

\section{The  degree of atypicality of a prime ideal.}

We can iterate the process described in Equations  \eqref{egg}, \eqref{elk} and Theorem \ref{cow}.  
Using the definitions of $ \A_{m,n}$  and $\B_{m,n}$ given near the start of section \ref{ai}, set $\A_i =  \gl  \A_{m-i,n-i},$ and 
$\B_i =    \B_{m-i,n-i},$ 
Then using the  $r+1$ pairs \[\{x_{m}, y_{n}\}, \ldots, \{x_{m-i}, y_{n-r}\},\]
 we obtain  surjective homomorphisms

\be \label{pug}  \A_0 \stackrel{\phi_1}{\lra} \A_{1}\stackrel{\phi_2}{\lra} \ldots \stackrel{\phi_r}{\lra} \A_{r}\ee
Set 
$\psi_i =\phi_i   \ldots \phi_2 \phi_1$, and in case S set
\[T_i= \prod^{m-i}_{k=1} \prod^{n-i}_{j=1}(x_k + y_j)\in     \A_{i}, \] $\B_i= \Z[\mathcal{X}_{m-i},\mathcal{Y}_{n-i}]^{W_i},$ and
 $\widehat{\B}_i= \Z[\mathcal{X}_{m-i},\mathcal{Y}_{n-i}]_{T_i}^{W_i},$   
where  ${W_i}$ is  a direct product of symmetric groups  ${W_i}= S_{m-i}\ti S_{n-i}$.   
\\ \\
In case L, set 
$\A_i =\Lambda_{m-i,n-i}$, 
$\B_i = \Lambda_{m-i} \ti \Lambda_{n-i}$, 
$$T_i=\prod_{j=1}^{n-i}\prod_{i=1}^{m-i}\left(1-\frac{y_{j}}{x_{i}}\right)
\in \A_{i},$$
 and $\widehat{\B}_i = (\Lambda_{m-i} \ti \Lambda_{n-i})_{T_i}$. 
Note that $T_i\B_i \subset    \A_i \subset \B_i$.
Repeating our earlier agruments  leads to a disjoint union of locally closed sets
\be \label{ram}
 \Spec   \A = \Spec\ B_T \;\cup\; \psi_1^{-1} (\Spec   
\widehat{\B}_{1}) \;\cup\;\ldots \;\cup\; \psi_{r-1}^{-1} (\Spec  \widehat{\B}_{r-1})\; \cup\; \psi_r^{-1} (\Spec  \A_{r}).\ee
Now suppose $p= \psi_r^{-1} (P)$, where $P \in   \Spec  \A_{r}$ and that $T_r \notin P$. In this case $P$ corresponds by localization to a prime ideal of $\widehat{\B}_r$, and 
 we say that $p$ has {\it  degree of atypicality} $r$.
 In the case that $P=M$ is a maximal ideal in $\A_r$, $M$ corresponds  to a maximal ideal $M_{T_r}$ of $\widehat{\B}_{r}$, and so 
$\cV(M)$ consists of a single  $W_{r}$  orbit of points in  $\Spec   \A_{r}.$
This shows that $p= \psi_r^{-1} (P)$ has degree of atyicality $r$ in the usual sense. 
Since  $\widehat{\B}_r$ is a Jacobson ring, we have the following. 
\bt    If $p$ is a prime ideal of $\cA$ with degree of atypicality $r$, then $p$ is an intersection of maximal ideals of $\cA$ each having degree of atypicality $r$. 
\et
\noi  Note also that if  $r= \min(m,n)$,  Equation \eqref{ram} gives  a natural stratification of $\Max \A_{m,n}$ according to the degree of atypicality.

\section{Applications to $Z(\fg)$.}
\subsection{The Image of the Harish-Chandra map.}
In this and the following subsection we apply our results to the center $Z(\fg)$ the enveloping algebra $U(\fg)$. Throughout we work over an algebraically closed field $\ttk$ of characteristic zero. 
For a basic classical simple Lie algebra the Harish-Chandra map, due to Gorelik and Kac,  gives an isomorphism $Z(\fg) \lra I(\fh)$, where $I(\fh)$ is a certain subalgebra of $S(\fh)^W$, for further details see \cite{M}  Theorem 13.1.1. 
In the next result we give a description of $I(\fh)$ when $\mathfrak{g}
= \fgl(m,n)$ or an orthosymplectic Lie superalgebra. Let $\epsilon_1, \ldots, \epsilon_m, \delta_1, \ldots, \delta_n$ be
the usual basis of $\mathfrak{h}^*$, and $h_1, \ldots, h_m,$ $ h'_1, \ldots,
h'_n$
 the dual basis for $\mathfrak{h}$, see \cite{M} Equations (2.2.4), (2.2.5) and (2.3.5).

If $\mathfrak{g} = \fgl(m,n)$ or $\mathfrak{g} =
\osp(2m+1,2n),$ $I(\fh)$ is easily expressed in terms of
supersymmetric polynomials. However if $\mathfrak{g} = \osp(2m,2n)$ the
situation  is more complicated.  In this case, set

$$J(\mathfrak{h}) = I_\ttk(h^2_1, \ldots , h^2_m; h'^{\; 2}_1, \ldots ,h'^{\; 2}_n), \quad \ T =\prod_{i,j}(h^2_i - (h'_j)^2)
, \quad \ \Phi = (h_1 \ldots h_m)T.$$ 
Also consider the group
   \[ W' = (\mathbb{Z}^m_2 \rtimes \mathcal{S}_m) \times (\mathbb{Z}^n_2
   \rtimes \mathcal{S}_n) \]
   where $\rtimes$ denotes a semidirect product.  There is an
   action of $W'$ on $S(\fh) = \ttk[h_1, \ldots, h_m, h'_1, \ldots,
   h'_n]$, where the symmetric groups $\mathcal{S}_m$ and $\mathcal{S}_n$ permute
   the $h_i, h'_i$ respectively, and $\mathbb{Z}^m_2,
   \mathbb{Z}^n_2$ change their signs.  The Weyl group $W$ of $\fg$ is a
   subgroup of index two in $W'.$ Namely $W$ consists of all elements of
   $W'$ which change an even number of signs of the $h_i$.
Then we have \cite{Se2}, \cite{M} Theorem 13.4.1.
 \bt
\label{ewe}  With the above notation  \bi \itema If
$\mathfrak{g} = \fgl(m,n)$ we have $I(\mathfrak{h})= I_\ttk(h_1, \ldots , h_m;
h_1', \ldots , h_n').$
\itemb If  $\mathfrak{g} = \osp(2m+1,2n)$
we have $I(\mathfrak{h}) =J(\mathfrak{h}).$
\itemc If $\mathfrak{g}  =  \osp(2m,2n),$
we have \ei 
 \be \label{kid} I(\mathfrak{h}) = J(\mathfrak{h}) +\Phi S(\fh)^{W'} .\ee  \et
\subsection{The case of $\osp(2m,2n)$.}
\noi The work on supersymmetric polynomials applies directly to $\Spec   I(\fh)$ when $\fg = \fgl(m,n)$ or $\osp(2m+1,2n)$.  Now assume that   $\mathfrak{g} = \osp(2m,2n)$. The first step is to observe that  the sum in Equation \eqref{kid} is  direct. Let 
$\gs:W'\lra \{\pm 1\}$ be the character with $\Ker \gs = W$.  If $R$ is a ring on which $W'$ acts,  such that $W$ acts trivially, we have $R = R_1 \op R_\gs$ where $R_1= R^{W'}$ and $$R_\gs = \{r \in R| wr = \gs(w)r \mbox{ for all  } w \in W'\}.$$ Note that 
$R_1$ is a subring of $R$ and $R_\gs$ is an ideal of $R$.
  If $C= S(\fh)^W$  we have $C_\gs =h_1 \ldots
h_m C_1$, and likewise for any  localization of $C$ with respect to a $W'$-invariant element. If $R= I(\mathfrak{h})$, we have $I(\mathfrak{h})_1 = J(\mathfrak{h})$ and 
$I(\mathfrak{h})_\gs=\Phi S(\fh)^{W'}$.
So the sum in \eqref{kid}
is just  $I(\mathfrak{h}) = I(\mathfrak{h})_1 \oplus I(\mathfrak{h})_\gs$.  

 Now we use the projection $\psi:I(\fh)\lra J(\fh)$ with kernel $I(\fh)_\gs$ to obtain an analog of Equation \eqref{gnu}.  
At this point the argument diverges slightly from the proof of Theorem \ref{cow}, so we give the details again.  
First observe  that for  any prime ideal $P$ of $I(\fh)$ with $T \in P$, we have $I(\fh)_\gs\subseteq P.$
Indeed   $ T$  is $W'$-invariant, and  
$T S(\fh)^{W'} \subset J(\mathfrak{h}) \subset I(\mathfrak{h})$. 
 Hence 
\be \label{asp} (I(\fh)_\gs)^2 =\Phi^2 S(\fh)^{W'} = T(Th_1^2 \ldots h_m^2) S(\fh)^{W'} \subseteq PI(\fh) = P. \ee
Furthermore

$$ J(\mathfrak{h})_{T}  = S(\fh)^{W'}_{T}\quad \mbox{ and } \quad \Phi S(\fh)^{W'}_{T}
=h_1 \ldots h_m S(\fh)^{W'}_{T},$$ so $ I(\mathfrak{h})_{T} = S(\fh)^{W}_{T}$.  Thus arguments similar to those leading to 
   \eqref{gnu} yield a disjoint union 
\[  \Spec   I(\fh) = \Spec   S(\fh)^W_{T}\cup \psi^{-1} (\Spec   J(\fh)).  \]
The analysis of maximal ideals is also slightly diffierent.  We use the  commutative diagram

\begin{xy} 
(0,20)*+{}="f";
(40,20)*+{ I(\fh)}="a";
 (100,20)*+{J(\fh)}="b"; (40,0)*+{ S(\fh)^W}="c"; 
(100,0)*+{ S(\fh)^{W'}}="d"; {\ar "a";"b"}?*!/_2mm/{\phi}; {\ar "a";"c"};
{\ar "b";"d"};
{\ar@{>} "c";"d"};?*!/_2mm/{\Phi};
 \end{xy}

\noi The map $\phi$ has kernel $I(\fh)_\gs$  and   $\Phi$ has kernel $S(\fh)_\gs^W$.   The two vertical maps are the obvious inclusions.

\bp If  $m\in \Max I(\fh)$, then $m{ S(\fh)^W}$ is a proper ideal of ${ S(\fh)^W}$. 
\ep
\bpf If $T\notin m$ this follows as in  Corollary \ref{pig}.  If  $T\in m$ then  
 $I(\fh)_\gs\subseteq m,$ as shown  in Equation \eqref{asp}, so $\phi(m)$ is 
a maximal ideal of $J(\fh)$. Since $J(\fh)$ is an algebra of supersymmetric polynomials, Proposition \ref{doe},  shows that  $L=\phi(m) S(\fh)^{W'}$ is a proper ideal of $ S(\fh)^{W'}$. 
 Because  $\Phi(mS(\fh)^W)=L,$ it follows that $mS(\fh)^W$ 
is a proper ideal of  $S(\fh)^W$.
\epf
\noi Using this result the weak and strong Nullstellensatz  for $I(\fh)$  may be deduced as before. In particular we have the following analog of Theorem \ref{hen}. 
\bt\label{rex} 
 If $m$ is a maximal ideal of 
$I(\mathfrak{h}))$, then there is a maximal ideal $M$ of  $S(\fh)^W$   such that 
$m = M \cap I(\mathfrak{h})$.
\et  
\noi Now Conjecture 13.5.1 from \cite{M} states that if $I(\fh)$ is as in Theorem \ref{ewe} and  $m$ is  any maximal ideal iof
$I(\fh)$, then there is a $\gl\in \fh^*$ such that 
\be \label{ox}
 m = \{ f \in I(\fh) | f(\gl) = 0\}.\ee
{\it Proof of the Conjecture.}
Since $I(\fh)$ is isomorphic to an algebra of supersymmetric polynomials in cases (a) and (b) of Theorem \ref{ewe}, we have  
$m = M \cap I(\fh)$ for some ideal of $S(\fh)^W$,   
by Theorem \ref{hen}.  By Theorem  \ref{rex} the same holds in case (c).   Now the maximal ideal $M$  in $S(\fh)^W$  corresponds to a $W$-orbit in $\fh^*,$ and Equation \eqref{ox} holds for any $\gl$ in this orbit.
\hfill $\Box$\\ \\

\section{Appendix: Characterizations of Laurent supersymmetric
 polynomials.}
We mention some alternative characterizations of Laurent supersymmetric
 polynomials with coefficients in  a field $K$.
Set  
\[ S=  K[x_1^{\pm1},\dots,x_m^{\pm1}|y_1^{\pm1},\dots,y_n^{\pm1}].\]
The direct product of symmetric groups $W = S_m\ti S_n$ acts on $S$ in the obvious way.
  Obviously if $f \in S$ is $W$-invariant, and
\be \label{ape} f(x_1 = t, y_1 = t)\; \mbox{is a polynomial which
is independent of}\; t,\ee
then $f(x_i = t, y_j = t)$ is
is independent of $t$ for all $i,j$.  We consider the condition (\ref{ape}) independently of $W$-invariance.
Let $\mathbb{T} $ be the algebraic torus $(K^*)^{m+n}$ and let $\{\epsilon_{1} , \ldots , \epsilon_{m}, \delta_{1} , \ldots , \delta_{n}\}$ be coordinates on $\mathbb{T}$.
Let $\ttT, \ttT_1 $ be the subtori of $\mathbb{T}$ defined by the equation $\ttT =\Ker \gep_1 \gd_1^{-1}, \ttT_1 =\Ker \gep_1 \cap \Ker\gd_1,$ and let $\ttT_2 = \{(t,1, \ldots,1| t,1, \ldots,1)|t\in K^*  \}.$   Then we have a direct product $\ttT =\ttT_1 \ttT_2$. 
Next set $R = K[x_2, \ldots, x_m| y_2, \ldots, y_n], x = x_1, y = y_1$, $z_+=(1 - \frac{x}{y}),$  and $z_-=(1 - \frac{y}{x}).$ 
  Note that $ (1 - z_+)^{-1}= (1-z_-)$ and hence 
$ z_- = z_+ (z_+-1)^{-1}$. 
It follows that  $S =  R[x^{\pm1}, y^{\pm1}]= R[x^{\pm1}, z_\pm],$   
$Sz_+=Sz_-$, and 

\be \label{pup} S=R[x^{\pm1},(1 - z_+)^{\pm1}].\ee
If $\gl =(x_1, x_2, \ldots, x_m| y_1, y_2, \ldots, y_n)$, $q\neq0\in K$ and $f \in S$, set 
$$f_q(\gl)=f(qx_1, x_2, \ldots, x_m| qy_1, y_2, \ldots, y_n).$$  
Then if  $\ga = (1, 0, \ldots, 0| 1, 0, \ldots, 0)$, we define
the {\it Laurent directional derivative} $D_\ga f$ in the direction of
$\exp(q\ga)$ by
\[ (D_\ga f)(\lambda) = \lim_{q \rightarrow 1}
\frac{f_q(\gl) - f(\gl)}{q-1}. \] This makes sense
since we only have to differentiate Laurent  polynomials. Note
that the directional derivative $D_\alpha f$ satisfies $D_\alpha f = x\partial f/\partial x + y\partial f/\partial y,$
where the partial derivatives vanish on $R.$ The result below is the Laurent analog of \cite{M} Lemma 12.1.1.

\begin{lemma}\label{cod}
Let $z =z_+$. For $f \in S$ the following conditions are equivalent
\begin{itemize} 
\itema $f \in R + Sz$
\itemb $f(x = y = t)$ is independent of $t\neq 0$ 
\itemc
For $\lambda \in \ttT$, and $t \in \ttT_2$, we have  $f(\lambda) = f(t\lambda)$.
\itemd $x\partial f/\partial x +y\partial f/\partial y \in (x - y)$
\end{itemize}
\end{lemma}
\noindent \bpf $(a) \Rightarrow (b)$ If $f \in R $ then $f(x = y =
t) = f$ is independent of $t$, while
if $f \in Sz$ then$f(x = y = t) = 0.$\\

$(a) \Rightarrow (d)$  This is similar to the proof of $(a)
\Rightarrow (b).$\\

$(c) \Leftrightarrow (d)$ 
This holds since  $D_\alpha f $  vanishes on
$\ttT $ if and only if $f_q({\lambda})$ is constant
for all $\lambda
\in \ttT .$\\

$(b) \Rightarrow (a)$ Using \eqref{pup}, we can write $f$ uniquely as a finite sum
\begin{equation} \label{fly}
f = \sum_{i \in \Z,j \in \Z}r_{i, j}x^i(1-{z})^j,
\end{equation}
 with $r_{i, j} \in R$ for
all $i,j.$ Then $f(x = y = t) = \sum_{i, j \in \Z}r_{i, j}t^i,$ so if
(b) holds then for all $i\neq 0$, we have  
\begin{equation} \label{emu} r_{i, 0} =-\sum_{j\neq 0} r_{i, j}.\ee Making this substitution in Equation  \eqref{fly} shows that $f \in R+ Sz+ Sz_{-} =R+ Sz.$\\

$(c) \Rightarrow (a)$ 
Given $\lambda \in \ttT $,
we can find $\lambda' \in \ttT_1 $ and $t' \in \ttT_2$ such that $\lambda =
t'\lambda'.$ Hence (c) is equivalent to $f(\lambda) =
f(t\lambda)$ for all $t \in \ttT_2$ and $\lambda \in \ttT_1 .$ 
Write  $f$ as in (\ref{fly}). Then for $\lambda \in \ttT_1 $ we have $ f(\lambda)  = \sum_{i \in \Z,j \in \Z}r_{i, j}(\gl)$ and $ f(t\lambda) = \sum_{i \in \Z,j \in \Z}r_{i, j}(\gl)t^i.$ 
Thus if (c) holds then again \eqref{emu} holds (with $r_{i, j}(\gl)$ replacing $r_{i, j}$) for all non-zero $i\in \Z,$ and (a) follows as before. \epf


\end{document}